\documentclass[11pt]{article}%
\usepackage{amssymb}
\usepackage{amsfonts}
\usepackage{amsmath}
\usepackage{graphicx}%
\providecommand{\U}[1]{\protect\rule{.1in}{.1in}}
\newtheorem{theorem}{Theorem}

\begin{document}

\title{Exact Markovian SIR and SIS epidemics on networks and an upper bound for the
epidemic threshold}
\author{P. Van Mieghem\thanks{ Faculty of Electrical Engineering, Mathematics and
Computer Science, P.O Box 5031, 2600 GA Delft, The Netherlands; \emph{email}:
P.F.A.VanMieghem@tudelft.nl }}
\date{Delft University of Technology}
\maketitle

\begin{abstract}
Exploiting the power of the expectation operator and indicator (or Bernoulli)
random variables, we present the exact governing equations for both the SIR
and SIS epidemic models on \emph{networks}. Although SIR and SIS are basic
epidemic models, deductions from their exact stochastic equations
\textbf{without} making approximations (such as the common mean-field
approximation) are scarce. An exact analytic solution of the governing
equations is highly unlikely to be found (for any network) due to the
appearing pair (and higher order) correlations. Nevertheless, the maximum
average fraction $y_{I}$ of infected nodes in both SIS and SIR can be written
as a quadratic form of the graph's Laplacian. Only for regular graphs, the
expression for the maximum of $y_{I}$ can be simplied to exhibit the explicit
dependence on the spectral radius. From our new Laplacian expression, we
deduce a general \textbf{upper} bound for the epidemic SIS threshold in any graph.

\end{abstract}

\section{Introduction}

Although the Susceptible-Infected-Removed (SIR) and the
Susceptible-Infected-Susceptible (SIS) model are basic corner-stones in
epidemics (see e.g.
\cite{Anderson_May,Daley,Newman_boek2010,Barrat_Bartelemy_Vespignani_CUPbook2008,Keeling_Rohani_2007,Diekmann_Heesterbeek_Britton_boek2012}%
), exact stochastic equations for SIR have, to the best of our knowledge, not
been published yet for an arbitrary network, while for SIS, we refer to
\cite{Simon_Taylor_Kiss_MathBiol2011} and
\cite{PVM_ToN_VirusSpread,PVM_EpsilonSIS_PRE2012}. A network is described by
an adjacency matrix $A$, with degree vector $D=\left(  d_{1},d_{2}%
,\ldots,d_{N}\right)  $ where $d_{k}$ is the degree of node $k$. For
simplicity, we assume an undirected network ($A=A^{T}$) that does not change
over time. In addition to the many applications ranging from cyber security
over information diffusion \cite{Vespignani_Science_2009} to biological
diseases \cite{Anderson_May,Keeling_Rohani_2007}, we explore these
(relatively) simple epidemic processes on graphs to understand the influence
of the topology of complex networks \cite{Albert_Barabasi_RevModPhys} on
properties of a dynamic process. First, we describe both the SIS and SIR model
on any network in a stochastic, Markovian setting and refer for non-Markovian
SIS epidemics to
\cite{PVM_nonMarkovianSIS_2013,PVM_nonMarkovianSIS_NIMFA_2013}.

In a SIS epidemic process, the viral state of a node $i$ at time $t$ is
specified by a Bernoulli random variable $X_{i}\left(  t\right)  \in\{S,I\}$:
$X_{i}\left(  t\right)  =S$ for a healthy, but susceptible node and
$X_{i}\left(  t\right)  =I$ for an infected node. A node $i$ at time $t$ can
be in one of the two states: \emph{infected}, with probability $v_{i}%
(t)=\Pr[X_{i}(t)=I]$ or \emph{healthy}, with probability $1-v_{i}(t)$, but
susceptible to the infection. We assume that the curing process per node $i$
is a Poisson process with rate $\delta$ and that the infection rate per link
is a Poisson process with rate $\beta$. Obviously, only when a node is
infected, it can infect its direct neighbors, that are still healthy. Both the
curing and infection Poisson process are independent. The effective infection
rate is defined by $\tau=\frac{\beta}{\delta}$. This is the general
continuous-time description of the simplest type of a SIS epidemic process on
a network.

In the SIR model, a node can be in one of the three states. When a node $j$ is
healthy, but susceptible to the virus, at time $t$, his state $Y_{j}=S$. A
node $j$ can be infected, $Y_{j}=I$, by its direct neighbors that are
infected. The infection is modelled by a Poisson process with rate $\beta$.
Finally, an infected node $j$ can be cured, after which it is removed from the
infection process, $Y_{j}=R$. The curing is modelled by a Poisson process with
rate $\delta$. All Poisson processes are independent. This formulation
describes a continuous-time SIR process on a graph.

There exist other formulations of the SIR process. For example, the
discrete-time counter part, in which a node is removed at the end of each
time-slot and infected neighbors can infect a susceptible node with
probability $p$, is termed a Reed-Frost process and is related to bond
percolation \cite{Newman_epidemics_PRE_2002}. Draief and Massouli\'{e}
\cite{Draief_Massoulie} show that a Reed-Frost process is related to the
growth of an Erd\H{o}s-R\'{e}nyi graph. The SIR process is also related to a
Markov discovery process on a graph (see \cite[p. 349-351]%
{PVM_PerformanceAnalysisCUP}). Newman \cite{Newman_epidemics_PRE_2002} has
presented a generating function approach for SIR, though implicitly assuming a
mean-field approximation. The above Markovian description of SIS\ and SIR,
based on independent Poisson processes, seems the most general one that still
allows us to write the general governing equations for any graph. Deviating
from a Markov process, by choosing other than the exponential interaction time
(for infection and/or curing, see
\cite{PVM_nonMarkovianSIS_2013,PVM_nonMarkovianSIS_NIMFA_2013}) or by
incorporating dependencies between the infection and curing process, will
complicate the analysis considerably. This argument provides the main
motivation to explore how far we can push the analysis to obtain physical insight.

\section{Governing equations}

In this paper, we analyse the SIR and SIS process rigorously and exploit the
power of the (linear) expectation operator $E\left[  .\right]  $ and the
indicator random variable $1_{x}$ (which equals one if the condition $x$ is
true, else it is zero) to remain closer to the physics of the epidemic
process. The SIR governing equation for the probability that a node $j$ is
infected reads%
\begin{equation}
\frac{d\Pr\left[  Y_{j}=I\right]  }{dt}=E\left[  -\delta1_{\left\{
Y_{j}=I\right\}  }+1_{\left\{  Y_{j}=S\right\}  }\beta\sum_{k=1}^{N}%
a_{kj}1_{\left\{  Y_{k}=I\right\}  }\right]  \label{dvg_prob_infected_node_j}%
\end{equation}
where the time-dependence of $Y_{j}\left(  t\right)  $ has been omitted for
simplicity. In words, the change in the probability that a node $j$ is
infected at time $t$ equals the expectation of (a) the rate $\beta$ times the
number of infected neighbors (specified by the adjacency matrix element
$a_{kj}$), given that node $j$ is susceptible minus (b) the rate $\delta$
given that the infected node is cured (and thereafter removed). Next, the
dynamic process that removes nodes satisfies%
\begin{equation}
\frac{d\Pr\left[  Y_{j}=R\right]  }{dt}=E\left[  \delta1_{\left\{
Y_{j}=I\right\}  }\right]  =\delta\Pr\left[  Y_{j}=I\right]
\label{dvg_prob_removed_node_j}%
\end{equation}
which says that the time-derivative of the probability that a node $j$ is
removed from the process equals the expectation of the rate $\delta$, given
that node $j$ is infected. Finally, a node is either healthy but susceptible,
infected, or cured (and removed); in other words, $1_{\left\{  Y_{j}%
=S\right\}  }+1_{\left\{  Y_{j}=I\right\}  }+1_{\left\{  Y_{j}=R\right\}  }=1$.

The first equation (\ref{dvg_prob_infected_node_j}) is complicating due to the
interaction with other infected nodes in the network, but
(\ref{dvg_prob_infected_node_j}) is of exactly the same form as the
corresponding SIS governing equation
\cite{PVM_secondorder_SISmeanfield_PRE2012},%
\[
\frac{d\Pr\left[  X_{j}=I\right]  }{dt}=E\left[  -\delta1_{\left\{
X_{j}=I\right\}  }+1_{\left\{  X_{j}=S\right\}  }\beta\sum_{k=1}^{N}%
a_{kj}1_{\left\{  X_{k}=I\right\}  }\right]
\]
However, in the SIS process, there are only two nodal states (or compartments)
possible so that $\ 1_{\left\{  X_{j}=S\right\}  }+1_{\left\{  X_{j}%
=I\right\}  }=1$, which leads to fewer equations than in the SIR process. We
proceed by rewriting equation (\ref{dvg_prob_infected_node_j}) using $E\left[
1_{\left\{  Y_{j}=S\right\}  \cap\left\{  Y_{k}=I\right\}  }\right]
=\Pr\left[  Y_{j}=S,Y_{k}=I\right]  $,%
\[
\frac{d\Pr\left[  Y_{j}=I\right]  }{dt}=-\delta\Pr\left[  Y_{j}=I\right]
+\beta\sum_{k=1}^{N}a_{kj}\Pr\left[  Y_{j}=S,Y_{k}=I\right]
\]
After invoking the law of total probability \cite[p. 27]%
{PVM_PerformanceAnalysisCUP},%
\[
\Pr\left[  Y_{k}=I\right]  =\Pr\left[  Y_{j}=S,Y_{k}=I\right]  +\Pr\left[
Y_{j}=I,Y_{k}=I\right]  +\Pr\left[  Y_{j}=R,Y_{k}=I\right]
\]
the SIR governing equation (\ref{dvg_prob_infected_node_j}) becomes%
\begin{equation}
\frac{d\!\Pr\!\left[  Y_{j}=I\right]  }{dt}=\beta\sum_{k=1}^{N}a_{kj}%
\!\Pr\!\left[  Y_{k}=I\right]  -\delta\Pr\!\left[  Y_{j}=I\right]  -\beta
\sum_{k=1}^{N}a_{kj}\!\left\{  \Pr\!\left[  Y_{j}=I,\!Y_{k}=I\right]
+\Pr\!\left[  Y_{j}=R,\!Y_{k}=I\right]  \right\}
\label{dvg_prob_infected_node_j_variant2}%
\end{equation}
The first two terms on the right-hand side in
(\ref{dvg_prob_infected_node_j_variant2}) describe the spread of the infection
from infected neighbors minus the nodal curing, while the third term excludes
infection spread to an infected or removed node $j$. This last term grows over
time, because (\ref{dvg_prob_removed_node_j}) illustrates that the probability
to become removed is non-decreasing over time. Relation
(\ref{dvg_prob_infected_node_j_variant2}) explains the bell-shape of
$\Pr\left[  Y_{j}\left(  t\right)  =I\right]  $ as a function of time $t$:
initially the third term is small and near to exponential growth arises from
the first and second term. As the number of removed nodes increases over time,
the third term counteracts the initial growth and forces its decline towards
extinction (for large $t$). The SIS differential equation corresponding to
(\ref{dvg_prob_infected_node_j_variant2}) is%
\begin{equation}
\frac{d\Pr\left[  X_{j}=I\right]  }{dt}=\beta\sum_{k=1}^{N}a_{kj}\Pr\left[
X_{k}=I\right]  -\delta\Pr\left[  X_{j}=I\right]  -\beta\sum_{k=1}^{N}%
a_{kj}\Pr\left[  X_{j}=I,X_{k}=I\right]  \label{dvg_prob_infected_node_j_SIS}%
\end{equation}

The governing equations (\ref{dvg_prob_infected_node_j_variant2}) and
(\ref{dvg_prob_infected_node_j_SIS}) lead to the following comparison: On the
same network under the same infection and curing rates and starting from one
infected node, the infection probability Pr$\left[  Y_{j}=I\right]  $ in SIR
epidemics is a lower bound for the infection probability $\Pr\left[
X_{j}=I\right]  $ in SIS epidemics. By starting the two processes on a same
network with the same initially infected node, the additional positive term
$\sum_{k=1}^{N}a_{kj}\!\Pr\!\left[  Y_{j}=R,\!Y_{k}=I\right]  $ in
(\ref{dvg_prob_infected_node_j_variant2}) shows that, at any time, Pr$\left[
Y_{j}=I\right]  \leq\Pr\left[  X_{j}=I\right]  $ for any node $j\in G$.
Physically, the removal process in SIR cannot increase the spread of infection
in the network with respect to SIS epidemics. As a consequence, the
$N$-intertwined mean-field approximation (NIMFA)
\cite{PVM_N_intertwined_Computing2011} upper bounds, besides SIS, also SIR epidemics.

Another interesting observation, also made in \cite{Parshani_PRL2010}, is that
the removal process in SIR epidemics prevents that a node can be infected
twice, which implies that the SIR infection process spreads over the network
as a growing discovery tree (without loops). Above the epidemic threshold,
most nodes are infected once (and discovered), while below the epidemic
threshold, the SIR infection tree dies out before infecting most nodes once.
Thus, in contrast to SIS epidemics, SIR infection travels from a node $i$ to a
node $j$ along a path, and not a walk. The tree spreading property of SIR
epidemics naturally maps SIR epidemics into a time-depending Bellman-Harris
branching process \cite{Harris_BP} on a network.

\section{Joint probabilities}

There are two ways to proceed from (\ref{dvg_prob_infected_node_j_variant2}):
either we deduce the governing equations for the two-pair probabilities as in
\cite{PVM_secondorder_SISmeanfield_PRE2012}, followed by higher order joint
probabilities until all $2^{N}$ SIS and $3^{N}$ SIR linear Markov equations
are established or we try to \textquotedblleft close\textquotedblright\ the
equations \cite[p. 653-654]{Newman_boek2010}, as coined in epidemiology. Here,
we propose a new method to compute all equations for higher order joint
probabilities. Indeed, interchanging the derivative and expectation operator
in (\ref{dvg_prob_infected_node_j}) yields
\begin{equation}
\frac{d1_{\left\{  X_{j}=I\right\}  }}{dt}=-\delta1_{\left\{  X_{j}=I\right\}
}+1_{\left\{  X_{j}=S\right\}  }\beta\sum_{k=1}^{N}a_{kj}1_{\left\{
X_{k}=I\right\}  } \label{definition_derivative_rv}%
\end{equation}
Strictly speaking, the derivative of an indicator does not exist, but we agree
to \emph{formally} define it by the random variable equation
(\ref{definition_derivative_rv}). Next, making the same reversal of operators,%
\[
\frac{d}{dt}E\left[
{\textstyle\prod\limits_{j=1}^{n}}
1_{\left\{  X_{j}=I\right\}  }\right]  \overset{^{\text{formally}}}{=}E\left[
\frac{d}{dt}%
{\textstyle\prod\limits_{j=1}^{n}}
1_{\left\{  X_{j}=I\right\}  }\right]  =E\left[  \sum_{m=1}^{n}%
{\textstyle\prod\limits_{j=1;j\neq m}^{n}}
1_{\left\{  X_{j}=I\right\}  }\frac{d1_{\left\{  X_{m}=I\right\}  }}%
{dt}\right]
\]
substituting (\ref{definition_derivative_rv}) and executing the $E\left[
.\right]  $ returns the correct result\footnote{The formal method can be made
mathematically rigorous (using the framework of stochastic differential
equations).},%
\begin{align*}
\frac{d}{dt}E\left[
{\textstyle\prod\limits_{j=1}^{n}}
1_{\left\{  X_{j}=I\right\}  }\right]   &  =-\delta nE\left[
{\textstyle\prod\limits_{j=1}^{n}}
1_{\left\{  X_{j}=I\right\}  }\right]  +\beta\sum_{m=1}^{n}\sum_{k=1}%
^{N}a_{km}E\left[  1_{\left\{  X_{k}=I\right\}  }%
{\textstyle\prod\limits_{j=1;j\neq m}^{n}}
1_{\left\{  X_{j}=I\right\}  }\right] \\
&  \hspace{0.5cm}-\beta\sum_{m=1}^{n}\sum_{k=1}^{N}a_{km}E\left[  1_{\left\{
X_{k}=I\right\}  }%
{\textstyle\prod\limits_{j=1}^{n}}
1_{\left\{  X_{j}=I\right\}  }\right]
\end{align*}
For each combination of $n$ out of $N$ states, such a differential equation
for the joint probability%
\[
E\left[
{\textstyle\prod\limits_{j=1}^{n}}
1_{\left\{  X_{j}=I\right\}  }\right]  =\Pr\left[  X_{1}=I,X_{2}%
=I,\ldots,X_{n}=I\right]
\]
can be written. The expectation in the last summation contains, except when
$\left(  1_{\left\{  X_{j}=I\right\}  }\right)  ^{2}=1_{\left\{
X_{j}=I\right\}  }$ occurs, a product of $n+1$ different random variables
$X_{j}$, for which a new differential equation is needed as outlined above. A
similar method applies for a product of different indicators, $%
{\textstyle\prod\limits_{j=1}^{n_{1}}}
1_{\left\{  Y_{j}=I\right\}  }%
{\textstyle\prod\limits_{j=n_{1}+1}^{n}}
1_{\left\{  Y_{j}=R\right\}  }$, where we define from
(\ref{dvg_prob_removed_node_j}) that $\frac{d1_{\left\{  Y_{j}=R\right\}  }%
}{dt}=\delta1_{\left\{  Y_{j}=I\right\}  }$. The analysis also shows that the
derivative of the $n$-th order joint probability includes joint probabilities
of order $n+1$, except if all nodes ($n=N$) are included and that an exact
description thus requires governing equations for all $1\leq n\leq N$ joint
probabilities, resulting in $2^{N}$ SIS and $3^{N}$ SIR linear Markov equations.

The most evident way of closure, which is an approximation method, is to
assume independence between nodes and states. For example, if we close the
first-order equations such as (\ref{dvg_prob_infected_node_j_variant2}) by
replacing $\Pr\left[  X_{j}=I,X_{k}=I\right]  $ by the product $f\left(
\Pr\left[  X_{j}=I\right]  \right)  g\left(  \Pr\left[  X_{k}=I\right]
\right)  $, where $f$ and $g$ are functions, we transform the set of linear
equations in first-order, $\Pr\left[  X_{m}=I\right]  $, and second-order,
$\Pr\left[  X_{j}=I,X_{k}=I\right]  $, variables to non-linear equations,
though with less variables (only first-order probabilities). This type of
approximation is also termed a mean-field approximation, that assumes
independence between the infection state of any two nodes.

\section{Properties deduced from first-order equations}

In the sequel, we continue to explore what can be deduced from the first-order
equations above without either higher-order deduction nor closure. We first
review a known result on the epidemic threshold for the SIS process that also
applies to the SIR process: The epidemic threshold of the SIR and
corresponding SIS process on any graph $G$ is lower bounded by
\begin{equation}
\tau_{c}\geq\frac{1}{\lambda_{1}}
\label{lower_bound_epidemic_threshold_SIS_SIR}%
\end{equation}
where $\lambda_{1}$ is the largest eigenvalue of the adjacency matrix $A$.
Directly from (\ref{dvg_prob_infected_node_j_variant2}) and
(\ref{dvg_prob_infected_node_j_SIS}), we deduce that%
\[
\frac{d\Pr\left[  Y_{j}\left(  t\right)  =I\right]  }{dt}\leq\beta\sum
_{k=1}^{N}a_{kj}\Pr\left[  Y_{k}\left(  t\right)  =I\right]  -\delta\Pr\left[
Y_{j}\left(  t\right)  =I\right]
\]
(and similarly for $\Pr\left[  X_{j}\left(  t\right)  =I\right]  $). The lower
bound (\ref{lower_bound_epidemic_threshold_SIS_SIR}) the follows by a similar
argument as in \cite{PVM_nonMarkovianSIS_2013}. The lower bound
(\ref{lower_bound_epidemic_threshold_SIS_SIR}) for the epidemic threshold also
holds for directed graphs. Since the SIR infection probability lower bounds
that of SIS in a same graph (with same initial conditions), $\tau_{c;SIS}%
\leq\tau_{c;SIR}$, which was earlier noted by Parshani \emph{et al.}
\cite{Parshani_PRL2010}.

For SIS epidemics, the lower bound
(\ref{lower_bound_epidemic_threshold_SIS_SIR}) was earlier proved in
\cite{PVM_ToN_VirusSpread}, though in a much less general and elegant form.
More importantly, the \textbf{lower} bound $\tau_{c}^{(1)}=\frac{1}%
{\lambda_{1}}$ appeared as the \textbf{exact} epidemic threshold in NIMFA,
where the superscript (1) in $\tau_{c}^{(1)}$ refers to the first order
mean-field approximation. We deem it important to underline the difference: in
the exact SIS and SIR model, the epidemic threshold $\tau_{c}$ is lower
bounded by $\tau_{c}^{(1)}=\frac{1}{\lambda_{1}}$, while in approximate
analyses (mean-field), the epidemic threshold is found to be equal to
$\tau_{c}^{(1)}=\frac{1}{\lambda_{1}}$. For some graphs (such as the complete
graph), the first order mean-field approximation $\tau_{c}^{(1)}$ is very
sharp, while for other graphs (such as the star), $\tau_{c}^{(1)}=\frac
{1}{\lambda_{1}}$ is less accurate \cite{PVM_MSIS_star_PRE2012}.

The \textbf{lower} bound $\tau_{c}^{(1)}=\frac{1}{\lambda_{1}}$ is of great
practical use: if the effective infection rate $\tau$ can be controlled such
that $\tau\leq\tau_{c}^{(1)}$ or the network can be designed to lower the
spectral radius $\lambda_{1}$ of a graph
\cite{PVM_decreasingspectralradius_PRE2011}, then the network is safeguarded
from long-term, massive infection. The lower bound
(\ref{lower_bound_epidemic_threshold_SIS_SIR}) cautions the widely cited
belief of a zero-epidemic threshold in scale-free networks
\cite{Boguna_Pastor_Vespignani_PRL_2003}: any finite network must have a
strictly positive epidemic threshold. Even when the mean-field epidemic
threshold $\tau_{c}^{(1)}\rightarrow0$ when $\lim_{N\rightarrow\infty}%
\lambda_{1}=\infty$, it may be possible, due to the \emph{lower} bound in
(\ref{lower_bound_epidemic_threshold_SIS_SIR}), that the exact threshold
$\tau_{c}>0$ is non-zero. An upper bound for $\frac{d\Pr\left[  X_{j}\left(
t\right)  =I\right]  }{dt}$ (and similarly for SIR) follows from the
H\"{o}lder inequality \cite[p. 90]{PVM_PerformanceAnalysisCUP} with $\frac
{1}{p}+\frac{1}{q}=1$ and $p>1$,%
\[
E\left[  1_{\left\{  X_{i}=I\right\}  }1_{\left\{  X_{k}=I\right\}  }\right]
\leq\left(  E\left[  1_{\left\{  X_{i}=I\right\}  }^{p}\right]  \right)
^{1/p}\left(  E\left[  1_{\left\{  X_{k}=I\right\}  }^{q}\right]  \right)
^{1/q}=\left(  \Pr\left[  X_{i}=I\right]  \right)  ^{1/p}\left(  \Pr\left[
X_{k}=I\right]  \right)  ^{1-1/p}%
\]
substituted into (\ref{dvg_prob_infected_node_j_SIS}) as%
\[
\frac{d\Pr\left[  X_{j}=I\right]  }{dt}\geq\beta\sum_{k=1}^{N}a_{kj}\Pr\left[
X_{k}=I\right]  -\delta\Pr\left[  X_{j}=I\right]  -\beta\left(  \Pr\left[
X_{i}=I\right]  \right)  ^{1/p}\sum_{k=1}^{N}a_{kj}\left(  \Pr\left[
X_{k}=I\right]  \right)  ^{1-1/p}%
\]
and the right-hand side can be maximized with respect to $p$. Unfortunately,
the steady-state solution of the above set of $N$ non-linear equations equals
$\Pr\left[  X_{j}=I\right]  =0$ for any node $j$ and any $p>1$. Recently,
Bogu$\widetilde{\text{n}}$a \emph{et al}.
\cite{Boguna_Castellano_Pastor_PRL2013} have proposed an approximate, coupling
type of argument to deduce an upper bound for the epidemic threshold. Although
their new method is ingenious and physically convincing, a proven upper bound
is still lacking. Below, we fill this gap by presenting a new and general
upper bound for the epidemic threshold $\tau_{c}$ in any network in Theorem
\ref{theorem_expression_epidemic_threshold} below.

By definition, the steady-state is attained for the time $t\rightarrow\infty$
at which the derivatives of the probabilities do not change anymore. If
$\frac{d\Pr\left[  Y_{j}=R\right]  }{dt}=0$ in (\ref{dvg_prob_removed_node_j})
for any node $j$, then $\Pr\left[  Y_{j}=I\right]  =0$ implying that there are
no infected nodes anymore in the network. In both SIS (due to the absorbing
state \cite{PVM_ToN_VirusSpread,PVM_EpsilonSIS_PRE2012}) and SIR epidemics,
the infectious disease eventually disappears from the network! Consequently,
the time-dependent (SIR) or metastable/quasi-stationary (SIS) behavior is
physically of interest. The final part expresses the exact prevalence in terms
of the graph's Laplacian $Q=\Delta-A$ (see e.g. \cite{PVM_graphspectra}) and
is proven in Appendix \ref{Proof_theorem_average_number_infected_nodes}:

\begin{theorem}
\label{theorem_average_number_infected_nodes} Denoting the (random) vector
$w_{I}=\left(  1_{\left\{  Y_{1}=I\right\}  },1_{\left\{  Y_{2}=I\right\}
},\ldots,1_{\left\{  Y_{N}=I\right\}  }\right)  $ and similarly for $w_{R}$,
the average number of infected nodes (or prevalence) satisfies for SIR
epidemics%
\begin{equation}
\frac{dy_{I}}{dt^{\ast}}=-y_{I}+\frac{\tau}{N}E\left[  w_{I}^{T}Qw_{I}%
-w_{I}^{T}Aw_{R}\right]  \label{dvgl_av_numb_infections_SIR}%
\end{equation}
while for SIS epidemics (denoted by a tilde)%
\begin{equation}
\frac{d\widetilde{y_{I}}}{dt^{\ast}}=-\widetilde{y_{I}}+\frac{\tau}{N}E\left[
\widetilde{w_{I}^{T}}Q\widetilde{w_{I}}\right]
\label{dvgl_av_numb_infections_SIS}%
\end{equation}
where $t^{\ast}=\delta t$ is the scaled time and $Q=\Delta-A$ is the Laplacian
of the graph with $\Delta=$ diag$\left(  d_{1},d_{2},\ldots,d_{N}\right)  $.
\end{theorem}

From (\ref{dvg_prob_removed_node_j}), we see that the average fraction of
removed nodes satisfies $\frac{dy_{R}}{dt^{\ast}}=y_{I}$. Apart from the
steady-state, also the maximum in (\ref{dvgl_av_numb_infections_SIR}) occurs
at $\frac{dy_{I}}{dt^{\ast}}=0$ and, at that value of time $t^{\ast}$, it
satisfies%
\begin{equation}
y_{I\max}=\frac{\tau}{N}E\left[  w_{I}^{T}Qw_{I}-w_{I}^{T}Aw_{R}\right]
\label{yImax_SIR}%
\end{equation}
illustrating that the corresponding $\widetilde{y}_{I\max}$ in SIS is larger
(because, in SIS, $w_{R}=0$ and $\widetilde{w}_{I}$ is not smaller on average
than $w_{I}$). In a regular graph, each node has degree $r$ and $Q=rI-A$ so
that (\ref{yImax_SIR}) simplifies to%
\[
y_{I\max}=\frac{\tau}{N}E\left[  rw_{I}^{T}w_{I}-w_{I}^{T}A(w_{I}%
+w_{R})\right]
\]
Since $w_{I}^{T}w_{I}=$ $\sum_{j=1}^{N}\left(  1_{\left\{  Y_{j}=I\right\}
}\right)  ^{2}=\sum_{j=1}^{N}1_{\left\{  Y_{j}=I\right\}  }=NZ_{I}$ and, thus
$y_{I}=NE\left[  w_{I}^{T}w_{I}\right]  $, we have%
\begin{equation}
y_{I\max}=\frac{\tau}{N}\frac{E\left[  w_{I}^{T}A(w_{I}+w_{R})\right]  }%
{r\tau-1} \label{yImax_SIR_regular_graph}%
\end{equation}
which illustrates (in agreement with
(\ref{lower_bound_epidemic_threshold_SIS_SIR}) because $\lambda_{1}=r$) that
$y_{I\max}=0$ when $\tau<\frac{1}{r}$ because $E\left[  w_{I}^{T}A(w_{I}%
+w_{R})\right]  \geq0$ and $y_{I}\geq0$. Only for regular graphs, the epidemic
threshold in both SIS and SIR epidemics appears directly from the exact
equation (\ref{yImax_SIR_regular_graph}). For special regular graphs such as
the complete graph, we can elaborate (\ref{yImax_SIR_regular_graph}) even
further. The natural extension from regular graphs to any graph is to bound
the degree vector as $d_{\min}u\leq D\leq d_{\max}u$ and
(\ref{dvgl_y_I_met_degree_vector}) becomes%
\[
\left(  \tau d_{\min}-1\right)  y_{I}-\frac{\tau}{N}E\left[  w_{I}^{T}%
A(w_{I}+w_{R})\right]  \leq\frac{dy_{I}}{dt^{\ast}}\leq\left(  \tau d_{\max
}-1\right)  y_{I}-\frac{\tau}{N}E\left[  w_{I}^{T}A(w_{I}+w_{R})\right]
\]
from which, for any graph, we find that%
\[
\frac{\tau}{N}\frac{E\left[  w_{I}^{T}A(w_{I}+w_{R})\right]  }{\tau d_{\max
}-1}\leq y_{I\max}\leq\frac{\tau}{N}\frac{E\left[  w_{I}^{T}A(w_{I}%
+w_{R})\right]  }{\tau d_{\min}-1}%
\]
illustrating, with (\ref{lower_bound_epidemic_threshold_SIS_SIR}), that the
epidemic threshold obeys $\frac{1}{d_{\max}}\leq\frac{1}{\lambda_{1}}\leq
\tau_{c}$. Since $E\left[  w_{I}^{T}A(w_{I}+w_{R})\right]  $ can still be zero
for $\tau>$ $\frac{1}{d_{\min}}$, we cannot conclude that $\tau_{c}\leq
\frac{1}{d_{\min}}$. In summary, a regular graph exhibits similar properties
as derived from mean-field or deterministic analyses. The larger the
heterogeneity in degree distribution as in most real-world networks
\cite{Albert_Barabasi_RevModPhys}, the larger we may expect that approximate
analyses deviate (see e.g. \cite{PVM_MSIS_star_PRE2012} for a star graph).

An upper bound for the SIS epidemic threshold, proven in Appendix
\ref{Proof_theorem_expression_epidemic_threshold}, is

\begin{theorem}
\label{theorem_expression_epidemic_threshold} Let $\varepsilon_{G}%
=\lim_{\widetilde{y_{I}}\downarrow0}\max_{\left(  k,l\right)  \in\mathcal{L}%
}\Pr\left[  X_{k}=I|X_{l}=I\right]  $, then the SIS epidemic threshold
$\tau_{c}$ in graph $G$ is upper bounded by%
\begin{equation}
\tau_{c}\leq\frac{1}{d_{\min}\left(  1-\varepsilon_{G}\right)  }
\label{upper_bound_epidemic_threshold}%
\end{equation}

\end{theorem}

The conditional probability $\varepsilon_{G}$ in Theorem
\ref{theorem_expression_epidemic_threshold} can be upper bounded by
$\varepsilon_{G}\leq\varepsilon_{K_{N}}$, because just at the onset of
infection ($\widetilde{y_{I}}\downarrow0$), the maximum conditional infection
probability $\varepsilon_{G}$ on a link $\left(  k,l\right)  $ in the graph
$G$ is largest in the complete graph. Exact computations on the complete graph
\cite{PVM_EpsilonSIS_PRE2012,PVM_MSIS_star_PRE2012} demonstrate that $\tau
_{c}=\frac{1}{N}\left(  1+\frac{c}{\sqrt{N}}+O\left(  N^{-1}\right)  \right)
$ for a constant $c$, implying that $\varepsilon_{K_{N}}=O\left(  \frac
{1}{\sqrt{N}}\right)  $ for large $N$. Hence, for large $N$, Theorem
\ref{theorem_expression_epidemic_threshold} leads to the upper bound%
\begin{equation}
\tau_{c}\leq\frac{1}{d_{\min}}\left(  1+O\left(  \frac{1}{\sqrt{N}}\right)
\right)  \label{upper_bound_epidemic_threshold_order_N}%
\end{equation}
for any graph\footnote{For large $N$, a lower bound for $\tau_{c}$ cannot be
of the form%
\[
\frac{1}{d_{\min}-x}%
\]
where $x$ is a fixed integer independently of $N$, because for the complete
graph $K_{N}$, $\frac{1}{d_{\min}-x}=\frac{1}{N-1-x}=\frac{1}{N}\left(
1+\frac{1+x}{N}+O\left(  N^{-2}\right)  \right)  $ which is smaller than the
exact threshold.}. Theorem \ref{theorem_expression_epidemic_threshold} (and
its proof) also emphasizes the role of the joint probability of infection at
end nodes of a same link, which laid at the basis of the pairwise
approximation \cite{Gleeson_PRL2011} and is considered as a significant
improvement over first-order mean-field approximations.

The upper bound (\ref{upper_bound_epidemic_threshold_order_N}) is sharp for
regular graphs, although (\ref{upper_bound_epidemic_threshold_order_N}) can be
large for realistic networks with broad (e.g. power law) degree distribution.
The general upper bound (\ref{upper_bound_epidemic_threshold}) and lower bound
(\ref{lower_bound_epidemic_threshold_SIS_SIR}) are, of course, less tight than
specific upper and lower bounds of particular classes of graphs, such as
regular trees, whose values are found in \cite[Table II]{Gleeson_PRL2011}
based on the work of Pemantle \cite{Pemantle_AP1992}, extended by Liggett
\cite{Liggett_AP1996}.

Finally, after tedious manipulations, the governing equation of the variance
of the fraction of infected nodes in SIS epidemics is%
\[
\frac{d\text{Var}\left[  \widetilde{Z}_{I}\right]  }{dt^{\ast}}=-2\text{Var}%
\left[  \widetilde{Z}_{I}\right]  +\frac{2\tau}{N}\left\{  E\left[
\widetilde{Z}_{I}\widetilde{w}_{I}^{T}Q\widetilde{w}_{I}\right]
-\widetilde{y}_{I}E\left[  \widetilde{w}_{I}^{T}Q\widetilde{w}_{I}\right]
\right\}  +\frac{1}{N}\left(  \widetilde{y}_{I}+\frac{\tau}{N}E\left[
\widetilde{w}_{I}^{T}Q\widetilde{w}_{I}\right]  \right)
\]
The variance is extremal when $\frac{d\text{Var}\left[  \widetilde{Z}%
_{I}\right]  }{dt^{\ast}}=0$, thus%
\begin{equation}
\left.  \text{Var}\left[  \widetilde{Z}_{I}\right]  \right\vert _{ex}%
=\frac{\tau}{N}\left\{  E\left[  \widetilde{Z}_{I}\widetilde{w}_{I}%
^{T}Q\widetilde{w}_{I}\right]  -\widetilde{y}_{I}E\left[  \widetilde{w}%
_{I}^{T}Q\widetilde{w}_{I}\right]  \right\}  +\frac{1}{2N}\left(
\widetilde{y}_{I}+\frac{\tau}{N}E\left[  \widetilde{w}_{I}^{T}Q\widetilde
{w}_{I}\right]  \right)  \label{Var[Z_I]_SIS_extremal}%
\end{equation}
The last term is never larger than $\frac{1}{N}$. If the fraction of infected
nodes $\widetilde{Z}_{I}$ and the sum over all links with precisely one end
infected, $\widetilde{w}_{I}^{T}Q\widetilde{w}_{I}=\sum_{l\in\mathcal{L}%
}\left(  1_{\left\{  X_{l^{+}}=I\right\}  }-1_{\left\{  X_{l^{-}}=I\right\}
}\right)  ^{2}$, were independent, then the maximum variance $\left.
\text{Var}\left[  \widetilde{Z}_{I}\right]  \right\vert _{ex}<\frac{1}{N}$
would be minimal. However, (\ref{dvgl_av_numb_infections_SIS}) shows that
$\widetilde{Z}_{I}$ and $\widetilde{w}_{I}^{T}Q\widetilde{w}_{I}$ are
dependent, implying that $\left.  \text{Var}\left[  \widetilde{Z}_{I}\right]
\right\vert _{ex}<1$ can be significant. For regular graphs,%
\[
\left.  \text{Var}\left[  \widetilde{Z}_{I}\right]  \right\vert _{ex}%
=\frac{E\left[  \widetilde{Z}_{I}\widetilde{w}_{I}^{T}A\widetilde{w}%
_{I}\right]  -\widetilde{y}_{I}E\left[  \widetilde{w}_{I}^{T}A\widetilde
{w}_{I}\right]  }{N\left(  r-\frac{1}{\tau}\right)  }+\frac{1}{2N}\left(
\frac{\frac{\tau}{N}E\left[  \widetilde{w}_{I}^{T}A\widetilde{w}_{I}\right]
-\widetilde{y}_{I}\left(  1+\tau r\right)  }{\tau r-1}\right)
\]
shows that the maximum variance occurs for $\tau$ around the epidemic
threshold $\tau_{c}\geq\frac{1}{r}$. The fact that the fraction of infected
nodes in SIS epidemics is found to vary most around the epidemic threshold,
where the process exhibits a phase transition (for large $N$), agrees with the
general physical theory of phase transitions \cite{Stanley_phase_transition}.

\section{Summary}

Based on the exact continuous-time, Markovian equations for SIS and SIR
epidemics, expressed in terms of Bernoulli random variables, we have proposed
a new method to deduce the differential equations for any joint probability.
Besides revisiting the known facts that the infection probability in SIS
epidemics always upper bounds that in SIR epidemics and that for both models,
the epidemic threshold is lower bounded by the inverse of the spectral radius,
we present a first order differential equation of the average SIS prevalence
over time containing the Laplacian of the graph, that elegantly expresses the
maximum average prevalence $y_{I\text{max}}$ in regular graphs in terms of the
spectral radius (or degree). From this new expression
(\ref{dvgl_av_numb_infections_SIS}), the SIS epidemic threshold in any graph
is upper bounded by (\ref{upper_bound_epidemic_threshold_order_N}), which
complements the result in \cite{Boguna_Castellano_Pastor_PRL2013}. Finally,
using our framework with Bernoulli random variables, the variance of the SIS
prevalence is computed and found to be maximal around the epidemic threshold.

\textbf{Acknowledgement.} The work is supported by EU project CONGAS (Grant
No. FP7-ICT-2011-8-317672). We are grateful to Eric Cator for the stimulating
discussions concerning Theorem \ref{theorem_expression_epidemic_threshold}.

{\footnotesize
\bibliographystyle{unsrt}
\bibliography{cac,MATH,misc,net,pvm,QTH,tel}
}

\newpage

\appendix{}

\section{Proof of the Theorems}

\subsection{Proof of Theorem \ref{theorem_average_number_infected_nodes}}

\label{Proof_theorem_average_number_infected_nodes}Summing
(\ref{dvg_prob_infected_node_j}) over all nodes $j$ yields%
\[
\frac{d}{dt}E\left[  \sum_{j=1}^{N}1_{\left\{  Y_{j}=I\right\}  }\right]
=E\left[  -\delta\sum_{j=1}^{N}1_{\left\{  Y_{j}=I\right\}  }+\beta\sum
_{k=1}^{N}1_{\left\{  Y_{k}=I\right\}  }\sum_{j=1}^{N}a_{kj}1_{\left\{
Y_{j}=S\right\}  }\right]
\]
Using $1_{\left\{  Y_{j}=S\right\}  }=1-1_{\left\{  Y_{j}=I\right\}
}-1_{\left\{  Y_{j}=R\right\}  }$, the last sum becomes%
\begin{align*}
\sum_{k=1}^{N}1_{\left\{  Y_{k}=I\right\}  }\sum_{j=1}^{N}a_{kj}1_{\left\{
Y_{j}=S\right\}  }  &  =\sum_{k=1}^{N}1_{\left\{  Y_{k}=I\right\}  }\left\{
\sum_{j=1}^{N}a_{kj}-\sum_{j=1}^{N}a_{kj}1_{\left\{  Y_{j}=I\right\}  }%
-\sum_{j=1}^{N}a_{kj}1_{\left\{  Y_{j}=R\right\}  }\right\} \\
&  =\sum_{k=1}^{N}d_{k}1_{\left\{  Y_{k}=I\right\}  }-\sum_{k=1}^{N}\sum
_{j=1}^{N}a_{kj}1_{\left\{  Y_{j}=I\right\}  }1_{\left\{  Y_{k}=I\right\}
}-\sum_{k=1}^{N}\sum_{j=1}^{N}a_{kj}1_{\left\{  Y_{k}=I\right\}  }1_{\left\{
Y_{j}=R\right\}  }\\
&  =D^{T}w_{I}-w_{I}^{T}Aw_{I}-w_{I}^{T}Aw_{R}%
\end{align*}
Further, denote by $Z_{I}=\frac{1}{N}\sum_{j=1}^{N}1_{\left\{  Y_{j}%
=I\right\}  }$ the fraction of infected nodes in the SIR process and by
$y_{I}=E\left[  Z_{I}\right]  $, then%
\[
N\frac{dy_{I}}{dt}=-N\delta y_{I}+\beta E\left[  D^{T}w_{I}-w_{I}^{T}%
Aw_{I}-w_{I}^{T}Aw_{R}\right]
\]
or, in terms of the effective infection rate $\tau=\frac{\beta}{\delta}$ in
units of $t^{\ast}=\delta t$,%
\begin{equation}
\frac{dy_{I}}{dt^{\ast}}=-y_{I}+\frac{\tau}{N}E\left[  D^{T}w_{I}-w_{I}%
^{T}Aw_{I}-w_{I}^{T}Aw_{R}\right]  \label{dvgl_y_I_met_degree_vector}%
\end{equation}
Using $D=\Delta u$, where $\Delta=$ diag$\left(  d_{1},d_{2},\ldots
,d_{N}\right)  $ and $u=(1,1,\ldots,1)$ is the all-one vector, we can rewrite%
\begin{align*}
D^{T}w_{I}-w_{I}^{T}Aw_{I}  &  =u^{T}\Delta w_{I}+w_{I}^{T}\Delta w_{I}%
-w_{I}^{T}\Delta w_{I}-w_{I}^{T}Aw_{I}\\
&  =\left(  u-w_{I}\right)  ^{T}\Delta w_{I}+w_{I}^{T}\left(  \Delta-A\right)
w_{I}%
\end{align*}
Since $1_{\left\{  Y_{j}=I\right\}  }1_{\left\{  Y_{j}=I\right\}
}=1_{\left\{  Y_{j}=I\right\}  }$,%
\[
\left(  u-w_{I}\right)  ^{T}\Delta w_{I}=\sum_{j=1}^{N}\left(  1-1_{\left\{
Y_{j}=I\right\}  }\right)  d_{j}1_{\left\{  Y_{j}=I\right\}  }=\sum_{j=1}%
^{N}\left(  1_{\left\{  Y_{j}=I\right\}  }-1_{\left\{  Y_{j}=I\right\}
}1_{\left\{  Y_{j}=I\right\}  }\right)  d_{j}=0
\]
Finally, introducing the Laplacian matrix $Q=\Delta-A$, we
arrive\footnote{Alternative expressions can be obtained using $u=w_{I}%
+w_{S}+w_{R}$ and $Qu=0$.} at (\ref{dvgl_av_numb_infections_SIR}). The SIS
variant (\ref{dvgl_av_numb_infections_SIS}) is similarly proved.\hfill
$\square\medskip$

\subsection{Proof of Theorem \ref{theorem_expression_epidemic_threshold}}

\label{Proof_theorem_expression_epidemic_threshold}From
(\ref{dvgl_av_numb_infections_SIS}) at $\frac{d\widetilde{y_{I}}}{dt^{\ast}%
}=0$, we find that%
\[
\tau^{-1}=\frac{E\left[  \widetilde{w_{I}^{T}}Q\widetilde{w_{I}}\right]
}{N\widetilde{y_{I}}}=\frac{E\left[  \widetilde{w_{I}^{T}}Q\widetilde{w_{I}%
}\right]  }{E\left[  \widetilde{w_{I}^{T}}u\right]  }%
\]
Introducing the basic Laplacian property $\widetilde{w_{I}^{T}}Q\widetilde
{w_{I}}=\sum_{l\in\mathcal{L}}\left(  1_{\left\{  X_{l^{+}}=I\right\}
}-1_{\left\{  X_{l^{-}}=I\right\}  }\right)  ^{2}$, where the link $l$ points
from node $l^{+}=i\rightarrow l^{-}=j$ and $\mathcal{L}$ is the set of links
of $G$, yields%
\begin{align*}
E\left[  \widetilde{w_{I}^{T}}Q\widetilde{w_{I}}\right]   &  =\sum
_{l\in\mathcal{L}}E\left[  \left(  1_{\left\{  X_{l^{+}}=I\right\}
}-1_{\left\{  X_{l^{-}}=I\right\}  }\right)  ^{2}\right]  =2\sum
_{l\in\mathcal{L}}E\left[  1_{\left\{  X_{l^{+}}=I\right\}  }-1_{\left\{
X_{l^{+}}=I\right\}  }1_{\left\{  X_{l^{-}}=I\right\}  }\right] \\
&  =2\sum_{l\in\mathcal{L}}E\left[  1_{\left\{  X_{l^{+}}=I\right\}  }\left(
1-1_{\left\{  X_{l^{-}}=I\right\}  }\right)  \right]  =2\sum_{l\in\mathcal{L}%
}E\left[  1_{\left\{  X_{l^{+}}=I\right\}  }1_{\left\{  X_{l^{-}}=S\right\}
}\right]
\end{align*}
Further, we can write%
\begin{align*}
2\sum_{l\in\mathcal{L}}E\left[  1_{\left\{  X_{l^{+}}=I\right\}  }1_{\left\{
X_{l^{-}}=S\right\}  }\right]   &  =2\sum_{l\in\mathcal{L}}\Pr\left[
X_{l^{+}}=I,X_{l^{-}}=S\right]  =\sum_{i=1}^{N}\sum_{j=1}^{N}a_{ij}\Pr\left[
X_{i}=I,X_{j}=S\right] \\
&  =\sum_{i=1}^{N}\Pr\left[  X_{i}=I\right]  \sum_{j=1}^{N}a_{ij}\Pr\left[
X_{j}=S|X_{i}=I\right]
\end{align*}
to obtain
\[
\tau^{-1}=\frac{\sum_{i=1}^{N}\Pr\left[  X_{i}=I\right]  \sum_{j=1}^{N}%
a_{ij}\Pr\left[  X_{j}=S|X_{i}=I\right]  }{\sum_{i=1}^{N}\Pr\left[
X_{i}=I\right]  }%
\]
The inequality \cite{Hardy_inequality}
\[
\min_{1\leq k\leq n}\frac{a_{k}}{q_{k}}\leq\frac{a_{1}+a_{2}+\cdots+a_{n}%
}{q_{1}+q_{2}+\cdots+q_{n}}\leq\max_{1\leq k\leq n}\frac{a_{k}}{q_{k}}%
\]
where $q_{1},q_{2},\ldots,q_{n}$ are positive real numbers and $a_{1}%
,a_{2},\ldots,a_{n}$ are real numbers leads to%
\[
\min_{1\leq i\leq N}\sum_{j=1}^{N}a_{ij}\Pr\left[  X_{j}=S|X_{i}=I\right]
\leq\tau^{-1}\leq\max_{1\leq i\leq N}\sum_{j=1}^{N}a_{ij}\Pr\left[
X_{j}=S|X_{i}=I\right]  \leq d_{\max}%
\]
Using the degree $d_{i}=\sum_{j=1}^{N}a_{ij}$, we proceed with the lower
bound,%
\begin{align*}
\tau^{-1}  &  \geq\min_{1\leq i\leq N}\sum_{j=1}^{N}a_{ij}\Pr\left[
X_{j}=S|X_{i}=I\right]  \geq\min_{1\leq i\leq N}\left(  \min_{\left(
k,l\right)  \in\mathcal{L}}\Pr\left[  X_{k}=S|X_{l}=I\right]  d_{i}\right) \\
&  =\min_{\left(  k,l\right)  \in\mathcal{L}}\Pr\left[  X_{k}=S|X_{l}%
=I\right]  d_{\min}%
\end{align*}
We define the epidemic threshold $\tau_{c}$ as that value of $\tau$ when the
prevalence (or order parameter) $\widetilde{y_{I}}=\frac{1}{N}\sum_{i=1}%
^{N}\Pr\left[  X_{i}=I\right]  $ approaches zero from above, denoted as
$\widetilde{y_{I}}\downarrow0$, so that%
\begin{equation}
\tau_{c}^{-1}=\lim_{\widetilde{y_{I}}\downarrow0}\frac{E\left[  \widetilde
{w_{I}^{T}}Q\widetilde{w_{I}}\right]  }{N\widetilde{y_{I}}} \label{def_tau_c}%
\end{equation}
and
\[
\tau_{c}^{-1}\geq d_{\min}\lim_{\widetilde{y_{I}}\downarrow0}\min_{\left(
k,l\right)  \in\mathcal{L}}\Pr\left[  X_{k}=S|X_{l}=I\right]
\]
The definition (\ref{def_tau_c}) of the epidemic threshold becomes
increasingly precise for large $N$. Finally, since $\Pr\left[  X_{k}%
=S|X_{l}=I\right]  =1-\Pr\left[  X_{k}=I|X_{l}=I\right]  $, we arrive at
(\ref{upper_bound_epidemic_threshold}). $\hfill\square\medskip$

\newpage

\section{The variance of $Z_{I}$ (in the SIS process)}

Recalling that the average fraction of infected nodes is $Z_{I}=\frac{1}%
{N}\sum_{j=1}^{N}1_{\left\{  X_{j}=I\right\}  }$ and omitting the tilde in the
notation (for SIS), then
\begin{align*}
E\left[  Z_{I}^{2}\right]   &  =\frac{1}{N^{2}}E\left[  \sum_{i=1}^{N}%
\sum_{j=1}^{N}1_{\left\{  X_{i}=I\right\}  \cap\left\{  X_{j}=I\right\}
}\right] \\
&  =\frac{1}{N^{2}}E\left[  \sum_{i=1}^{N}\sum_{j=1;j\neq i}^{N}1_{\left\{
X_{i}=I\right\}  \cap\left\{  X_{j}=I\right\}  }+\sum_{i=1}^{N}1_{\left\{
X_{i}=I\right\}  \cap\left\{  X_{i}=I\right\}  }\right] \\
&  =\frac{s_{I}}{N^{2}}+\frac{y_{I}}{N}%
\end{align*}
where%
\[
s_{I}=E\left[  \sum_{i=1}^{N}\sum_{j=1;j\neq i}^{N}1_{\left\{  X_{i}%
=I\right\}  \cap\left\{  X_{j}=I\right\}  }\right]
\]
First, using (\ref{definition_derivative_rv}), we have, for $i\neq j$,
\begin{align}
\frac{d}{dt}E\left[  1_{\left\{  X_{i}=I\right\}  \cap\left\{  X_{j}%
=I\right\}  }\right]   &  =E\left[  1_{\left\{  X_{i}=I\right\}  }%
\frac{d1_{\left\{  X_{j}=I\right\}  }}{dt}+1_{\left\{  X_{j}=I\right\}  }%
\frac{d1_{\left\{  X_{i}=I\right\}  }}{dt}\right] \nonumber\\
&  =E\left[  -2\delta1_{\left\{  X_{i}=I\right\}  \cap\left\{  X_{j}%
=I\right\}  }+\beta1_{\left\{  X_{i}=I\right\}  \cap\left\{  X_{j}=S\right\}
}\sum_{k=1}^{N}a_{kj}1_{\left\{  X_{k}=I\right\}  }\right. \nonumber\\
&  \hspace{0.5cm}\left.  +\beta1_{\left\{  X_{j}=I\right\}  \cap\left\{
X_{i}=S\right\}  }\sum_{k=1}^{N}a_{ki}1_{\left\{  X_{k}=I\right\}  }\right]
\label{dvg_joint_prob}%
\end{align}
Summing over all $i$ and $j\neq i$ yields, in time units of $t^{\ast}=\delta
t$,%
\[
\frac{ds_{I}}{dt^{\ast}}=-2s_{I}+\tau E\left[  \sum_{i=1}^{N}\sum_{j=1;j\neq
i}^{N}1_{\left\{  X_{i}=I\right\}  \cap\left\{  X_{j}=S\right\}  }\sum
_{k=1}^{N}a_{kj}1_{\left\{  X_{k}=I\right\}  }+\sum_{i=1}^{N}\sum_{j=1;j\neq
i}^{N}1_{\left\{  X_{j}=I\right\}  \cap\left\{  X_{i}=S\right\}  }\sum
_{k=1}^{N}a_{ki}1_{\left\{  X_{k}=I\right\}  }\right]
\]
Using $1_{\left\{  X_{j}=S\right\}  }=1-1_{\left\{  X_{j}=I\right\}  }$, we
have that%
\begin{align*}
R  &  =\sum_{i=1}^{N}\sum_{j=1;j\neq i}^{N}1_{\left\{  X_{i}=I\right\}
\cap\left\{  X_{j}=S\right\}  }\sum_{k=1}^{N}a_{kj}1_{\left\{  X_{k}%
=I\right\}  }+\sum_{i=1}^{N}\sum_{j=1;j\neq i}^{N}1_{\left\{  X_{j}=I\right\}
\cap\left\{  X_{i}=S\right\}  }\sum_{k=1}^{N}a_{ki}1_{\left\{  X_{k}%
=I\right\}  }\\
&  =\sum_{i=1}^{N}\sum_{j=1;j\neq i}^{N}\sum_{k=1}^{N}a_{kj}1_{\left\{
X_{i}=I\right\}  }\left(  1-1_{\left\{  X_{j}=I\right\}  }\right)  1_{\left\{
X_{k}=I\right\}  }+\sum_{i=1}^{N}\sum_{j=1;j\neq i}^{N}\sum_{k=1}^{N}%
a_{ki}1_{\left\{  X_{j}=I\right\}  }\left(  1-1_{\left\{  X_{i}=I\right\}
}\right)  1_{\left\{  X_{k}=I\right\}  }\\
&  =\sum_{i=1}^{N}\sum_{j=1;j\neq i}^{N}\sum_{k=1}^{N}a_{kj}1_{\left\{
X_{i}=I\right\}  }1_{\left\{  X_{k}=I\right\}  }-\sum_{i=1}^{N}\sum_{j=1;j\neq
i}^{N}\sum_{k=1}^{N}a_{kj}1_{\left\{  X_{i}=I\right\}  }1_{\left\{
X_{j}=I\right\}  }1_{\left\{  X_{k}=I\right\}  }\\
&  +\sum_{i=1}^{N}\sum_{j=1;j\neq i}^{N}\sum_{k=1}^{N}a_{ki}1_{\left\{
X_{j}=I\right\}  }1_{\left\{  X_{k}=I\right\}  }-\sum_{i=1}^{N}\sum_{j=1;j\neq
i}^{N}\sum_{k=1}^{N}a_{ki}1_{\left\{  X_{j}=I\right\}  }1_{\left\{
X_{i}=I\right\}  }1_{\left\{  X_{k}=I\right\}  }%
\end{align*}
Further,%
\begin{align*}
R  &  =\sum_{i=1}^{N}1_{\left\{  X_{i}=I\right\}  }\sum_{k=1}^{N}1_{\left\{
X_{k}=I\right\}  }\sum_{j=1;j\neq i}^{N}a_{kj}-\sum_{i=1}^{N}1_{\left\{
X_{i}=I\right\}  }\sum_{k=1}^{N}1_{\left\{  X_{k}=I\right\}  }\sum_{j=1;j\neq
i}^{N}a_{kj}1_{\left\{  X_{j}=I\right\}  }\\
&  +\sum_{k=1}^{N}1_{\left\{  X_{k}=I\right\}  }\sum_{i=1}^{N}a_{ki}%
\sum_{j=1;j\neq i}^{N}1_{\left\{  X_{j}=I\right\}  }-\sum_{i=1}^{N}\sum
_{k=1}^{N}a_{ki}1_{\left\{  X_{i}=I\right\}  }1_{\left\{  X_{k}=I\right\}
}\sum_{j=1;j\neq i}^{N}1_{\left\{  X_{j}=I\right\}  }\\
&  =\sum_{i=1}^{N}1_{\left\{  X_{i}=I\right\}  }\sum_{k=1}^{N}1_{\left\{
X_{k}=I\right\}  }\left(  \sum_{j=1}^{N}a_{kj}-a_{ki}\right)  -\sum_{i=1}%
^{N}1_{\left\{  X_{i}=I\right\}  }\sum_{k=1}^{N}1_{\left\{  X_{k}=I\right\}
}\left(  \sum_{j=1}^{N}a_{kj}1_{\left\{  X_{j}=I\right\}  }-a_{ki}1_{\left\{
X_{i}=I\right\}  }\right) \\
&  +\sum_{k=1}^{N}1_{\left\{  X_{k}=I\right\}  }\sum_{i=1}^{N}a_{ki}\left(
\sum_{j=1}^{N}1_{\left\{  X_{j}=I\right\}  }-1_{\left\{  X_{i}=I\right\}
}\right)  -\sum_{i=1}^{N}\sum_{k=1}^{N}a_{ki}1_{\left\{  X_{i}=I\right\}
}1_{\left\{  X_{k}=I\right\}  }\left(  \sum_{j=1}^{N}1_{\left\{
X_{j}=I\right\}  }-1_{\left\{  X_{i}=I\right\}  }\right) \\
&  =\sum_{i=1}^{N}1_{\left\{  X_{i}=I\right\}  }\sum_{k=1}^{N}1_{\left\{
X_{k}=I\right\}  }d_{k}-\sum_{i=1}^{N}\sum_{k=1}^{N}1_{\left\{  X_{k}%
=I\right\}  }a_{ki}1_{\left\{  X_{i}=I\right\}  }-\sum_{i=1}^{N}1_{\left\{
X_{i}=I\right\}  }\sum_{k=1}^{N}1_{\left\{  X_{k}=I\right\}  }\sum_{j=1}%
^{N}a_{kj}1_{\left\{  X_{j}=I\right\}  }\\
&  +\sum_{i=1}^{N}\sum_{k=1}^{N}1_{\left\{  X_{k}=I\right\}  }a_{ki}%
1_{\left\{  X_{i}=I\right\}  }+\sum_{k=1}^{N}1_{\left\{  X_{k}=I\right\}
}d_{k}\left(  \sum_{j=1}^{N}1_{\left\{  X_{j}=I\right\}  }\right)  -\sum
_{k=1}^{N}\sum_{i=1}^{N}1_{\left\{  X_{k}=I\right\}  }a_{ki}1_{\left\{
X_{i}=I\right\}  }\\
&  -\sum_{i=1}^{N}\sum_{k=1}^{N}a_{ki}1_{\left\{  X_{i}=I\right\}
}1_{\left\{  X_{k}=I\right\}  }\sum_{j=1}^{N}1_{\left\{  X_{j}=I\right\}
}+\sum_{i=1}^{N}\sum_{k=1}^{N}a_{ki}1_{\left\{  X_{i}=I\right\}  }1_{\left\{
X_{k}=I\right\}  }%
\end{align*}
Hence,%
\begin{align*}
R  &  =2\sum_{i=1}^{N}1_{\left\{  X_{i}=I\right\}  }\sum_{k=1}^{N}1_{\left\{
X_{k}=I\right\}  }d_{k}-\sum_{i=1}^{N}1_{\left\{  X_{i}=I\right\}  }\sum
_{k=1}^{N}\sum_{j=1}^{N}1_{\left\{  X_{k}=I\right\}  }a_{kj}1_{\left\{
X_{j}=I\right\}  }\\
&  -\sum_{j=1}^{N}1_{\left\{  X_{j}=I\right\}  }\sum_{i=1}^{N}\sum_{k=1}%
^{N}1_{\left\{  X_{k}=I\right\}  }a_{ki}1_{\left\{  X_{i}=I\right\}  }%
\end{align*}
and in vector form,%
\[
R=2NZ_{I}\left(  D^{T}w_{I}-w_{I}^{T}Aw_{I}\right)
\]
Combining all parts, with
\[
D^{T}w_{I}-w_{I}^{T}Aw_{I}=w_{I}^{T}Qw_{I}%
\]
as shown in the proof of Theorem \ref{theorem_average_number_infected_nodes},
we have%
\[
\frac{ds_{I}}{dt}=-2\delta s_{I}+2N\beta E\left[  Z_{I}w_{I}^{T}Qw_{I}\right]
\]
so that%
\begin{align*}
\frac{dE\left[  Z_{I}^{2}\right]  }{dt}  &  =\frac{d}{dt}\left(  \frac{s_{I}%
}{N^{2}}+\frac{y_{I}}{N}\right)  =\frac{1}{N^{2}}\frac{ds_{I}}{dt}+\frac{1}%
{N}\frac{dy_{I}}{dt}\\
&  =\frac{1}{N^{2}}\left(  -2\delta s_{I}+2N\beta E\left[  Z_{I}w_{I}%
^{T}Qw_{I}\right]  \right)  +\frac{1}{N}\left(  -\delta y_{I}+\frac{\beta}%
{N}E\left[  w_{I}^{T}Qw_{I}\right]  \right) \\
&  =\frac{-\delta\left(  2s_{I}+Ny_{I}\right)  }{N^{2}}+\frac{2\beta}%
{N}E\left[  Z_{I}w_{I}^{T}Qw_{I}\right]  +\frac{\beta}{N^{2}}E\left[
w_{I}^{T}Qw_{I}\right]
\end{align*}
Finally, Var$\left[  Z_{I}\right]  =E\left[  Z_{I}^{2}\right]  -y_{I}^{2}$,
from which%
\begin{align*}
\frac{d\text{Var}\left[  Z_{I}\right]  }{dt}  &  =\frac{dE\left[  Z_{I}%
^{2}\right]  }{dt}-2y_{I}\frac{dy_{I}}{dt}\\
&  =\frac{-\delta\left(  2s_{I}+Ny_{I}\right)  }{N^{2}}+\frac{2\beta}%
{N}E\left[  Z_{I}w_{I}^{T}Qw_{I}\right]  +\frac{\beta}{N^{2}}E\left[
w_{I}^{T}Qw_{I}\right]  -2y_{I}\left(  -\delta y_{I}+\frac{\beta}{N}E\left[
w_{I}^{T}Qw_{I}\right]  \right) \\
&  =\delta\frac{2N^{2}y_{I}^{2}-\left(  2s_{I}+Ny_{I}\right)  }{N^{2}}%
+\frac{2\beta}{N}\left\{  E\left[  Z_{I}w_{I}^{T}Qw_{I}\right]  -y_{I}E\left[
w_{I}^{T}Qw_{I}\right]  \right\}  +\frac{\beta}{N^{2}}E\left[  w_{I}^{T}%
Qw_{I}\right]
\end{align*}
Now,%
\begin{align*}
\frac{2N^{2}y_{I}^{2}-\left(  2s_{I}+Ny_{I}\right)  }{N^{2}}  &  =\frac
{2N^{2}y_{I}^{2}-2N^{2}E\left[  Z_{I}^{2}\right]  +Ny_{I}}{N^{2}}\\
&  =-2\text{Var}\left[  Z_{I}\right]  +\frac{y_{I}}{N}%
\end{align*}
Thus,%
\[
\frac{d\text{Var}\left[  Z_{I}\right]  }{dt^{\ast}}=-2\text{Var}\left[
Z_{I}\right]  +\frac{2\tau}{N}\left\{  E\left[  Z_{I}w_{I}^{T}Qw_{I}\right]
-y_{I}E\left[  w_{I}^{T}Qw_{I}\right]  \right\}  +\frac{1}{N}\left(
y_{I}+\frac{\tau}{N}E\left[  w_{I}^{T}Qw_{I}\right]  \right)
\]

\end{document}